\documentclass{amsart}

\usepackage{amsthm}
\usepackage{amsfonts}

\newcommand{\C}{{\mathbb{C}}}
\newcommand{\Pp}{{\mathbb{P}}}

\newcommand{\bo}{{\mathfrak{b}}}
\newcommand{\gp}{{\mathfrak{g}}}

\newcommand{\tor}{{\mathfrak{t}}}
\newcommand{\He}{{\mathcal{H}}}
\newcommand{\Cy}{{Y}}
\newcommand{\X}[1]{{\mathcal{X}}_{#1}}
\newcommand{\St}{{\textup{Stab}}}

\newtheorem{theorem}{Theorem}[section]
\newtheorem{corollary}[theorem]{Corollary}
\newtheorem{lemma}[theorem]{Lemma}
\newtheorem{proposition}[theorem]{Proposition}

\newtheorem{definition}[theorem]{Definition}

\newtheorem{question}[theorem]{Question}

\title{Hessenberg varieties are not pure dimensional}
\author{Julianna S. Tymoczko}

\begin{document}
\begin{abstract}
We study a family of subvarieties of the flag variety defined by certain 
linear conditions, called Hessenberg varieties.  We compare them to Schubert varieties.
We prove that some Schubert varieties can be realized as Hessenberg varieties and
vice versa.  Our proof explicitly identifies these Schubert varieties by their
permutation and computes their dimension.

We use this to answer an open question by proving 
that Hessenberg varieties are not always pure dimensional.  We give
examples that neither semisimple nor nilpotent Hessenberg varieties need
be pure; the latter are connected, non-pure-dimensional
Hessenberg varieties.  Our methods require us to generalize the definition of Hessenberg
varieties.
\end{abstract}

\maketitle

\begin{center} {\em Dedicated to Robert MacPherson.}  \end{center}

\section{Introduction: Background and notation}

A flag is a nested collection of vector spaces $V_1 \subseteq V_2 \subseteq \cdots
\subseteq V_n = \C^n$, where each $V_i$ is $i$-dimensional.  The full flag variety
is the complex algebraic variety consisting of all flags; it is smooth and compact.

This paper studies two families of subvarieties of the full flag variety: Hessenberg varieties
and Schubert varieties.  The first family is defined using two parameters: a linear operator $X: \C^n
\rightarrow \C^n$ and a nondecreasing function $h: \{1,2,\ldots, n\} \longrightarrow
\{1,2,\ldots,n\}$.  We call $h$ a Hessenberg function.  The Hessenberg 
variety associated to $X$ and $h$ is denoted $\He(X,h)$ and defined by
\[\He(X,h) = \{{\textup{Flags }}: XV_i \subseteq V_{h(i)} \textup{ for all }i\}.\]
(This generalizes the original definition of
\cite{dMPS}, as in Sections
\ref{hess vars} and \ref{generalizations}.)  

 For instance, if $X$ is arbitrary and $h$ has $h(i)=n$ 
for all $i$, then $\He(X,h)$
is the full flag variety.   More interesting are the Springer fibers, namely the Hessenberg varieties 
such that $X$ is nilpotent and $h(i)=i$ for each $i$.  Springer fibers are used to 
construct geometric representations of the symmetric group (\cite{CG} gives a survey).  
W.\ Borho and R.\ MacPherson generalize Springer representations to
a class of Hessenberg varieties that blend these two examples:
 $h$ is a {\em parabolic} function, defined in Section \ref{generalizations}, and 
 $X$ is a nilpotent matrix whose Jordan blocks are subordinate to $h$ (see \cite{BM}).  
 More Hessenberg varieties are in Section \ref{hess vars}.

This paper answers an open question about Hessenberg varieties: are they all
pure dimensional?  The pure-dimensionality
of Springer fibers is significant for Springer representations, which arise from permutation actions on top-dimensional cohomology.   Until now, the answer was yes in all known cases.  

We
show two ways in which Hessenberg varieties can fail to be pure dimensional.
In Section \ref{semisimple impurity}, we give an example in which $X$ is a semisimple 
operator and $\He(X,h)$ is a disjoint union of smooth subvarieties of $G/B$ of different
dimensions.  One case of this example came up in calculations that R.\ MacPherson 
and I performed while researching \cite{MT}.
In Section \ref{highweight}, we show that $\He(X,h)$ need not be
pure dimensional even when $X$ is nilpotent.  Section \ref{highweight} 
gives a family of examples 
that are connected but (in general) reducible
Hessenberg varieties whose components have different dimensions.

To prove that nilpotent Hessenberg varieties are not always pure, we
use Schubert varieties.  
Every invertible matrix $g$ gives a flag $[g]$ whose $i$-dimensional 
subspace is spanned by the first $i$ columns of $g$.  For each permutation $w$, the
Schubert variety
$\Cy_w$ is the closure of the set $\{[bw]: b \textup{ is upper-triangular}\}$.  
Schubert varieties are important 
because they form a basis for the cohomology of the flag variety.
Their geometry is a subject of intense scrutiny and is related 
to the combinatorics of the symmetric group.
For instance, whether $\Cy_w$ is singular is determined by 
substrings of  $w$ 
\cite[Chapters 5 and 8]{BL}.  

We show
that certain Schubert varieties can be realized as Hessenberg varieties, and conversely that 
some Hessenberg varieties are unions of Schubert varieties.  To construct
these Schubert varieties, we take $X$ to be the highest weight vector, namely
$X = E_{1n}$.  Section \ref{highweight} describes these Hessenberg varieties
in terms of their Schubert-variety components.

Most of this paper treats full flags in $GL_n(\C)$.  
Section \ref{generalizations} discusses how to generalize these results to other Lie types.
Section \ref{questions} contains open questions about Hessenberg varieties, including
the question of whether every Schubert variety can be realized as a Hessenberg
variety.

The author thanks Konstanze Rietsch and John Stembridge
for inspiring and helpful conversations.

\subsection{Descriptions of the flag variety}

This section is primarily an exposition of three classical ways to describe the flag variety, one
geometric, one algebraic, and one combinatorial, all three of which will be used in this paper.
This section also includes small
lemmas needed elsewhere.  Our motivation when selecting these proofs was diversity
of approach.

\subsubsection{Geometric description of the flag variety}

Our initial definition was a geometric characterization of the variety of full flags in $\C^n$.  
We denote the flag $V_1 \subseteq \cdots \subseteq V_n$ by $V_{ \bullet}$.  

Throughout this paper, we use a fixed 
basis $e_1$, $\ldots$, $e_n$ for $\C^n$.  Each flag can be written explicitly in terms
of this basis.
%Bott-Samelson resolution could be introduced here, but we won't use it in this paper.
%If it is described, describe length of permutation in terms of Bott-Samelson resolutions

\subsubsection{Algebraic description of the flag variety}

The flag $V_{\bullet}$ can be realized (non-uniquely) as an invertible matrix $g$ using the 
rule that the first $i$ columns of $g$ span $V_i$.  In this case $V_{\bullet}$ is also denoted $[g]$.  Let $B$ denote the group of invertible
upper-triangular matrices.  The flag variety is the
quotient $GL_n/B$.

The group $GL_n$ acts on the flag variety by the rule 
that if $h$ is in $GL_n$ and $[g]$ is in $GL_n/B$ then $h \cdot [g] = [hg]$.  
When this action is restricted to the upper-triangular matrices $B$, it partitions the flag
variety into $B$-orbits whose closures are the Schubert varieties $\Cy_w$.  
%When this action is restricted
%to the diagonal matrices $T$, it gives a torus action with finitely many fixed points and
%finitely many one-dimensional orbits.  
%\begin{todo} check: fixed points and one-orbits necessary? \end{todo}
%Our next description of the flag variety helps
%us work with these orbits.

\subsubsection{Combinatorial description of the flag variety}

The permutation matrices index Schubert varieties and, as flags, are contained in 
Schubert varieties.  We use $w$ to refer both to the permutation matrix and to the permutation
on the set $\{1,2,\ldots,n\}$ defined by $we_i = e_{w(i)}$.  We denote
transpositions by $s_{ij}$ and denote arbitrary permutations by $w$ or $v$.  

For each $w$, the {\em Schubert cell}
$[Bw]$ is the interior of the Schubert variety $\Cy_w$.  It can be described explicitly 
using the following subgroup of $B$.

\begin{definition}
Fix a permutation $w$ and let $U_w$ be the subgroup of $B$ 
defined by either one of the following equivalent conditions:
\begin{enumerate}
\item $U_w$ is the maximal subgroup of $B$ such that $w^{-1}U_w w$ is lower-triangular
	with ones along the diagonal.
\item $U_w$ consists of all matrices in $B$ with ones along the diagonal and whose
	$(i,j)$ entry is zero for each pair $i<j$ with $w^{-1}(i)<w^{-1}(j)$.
\end{enumerate}
\end{definition}

The next proposition follows from \cite[Sections 28.3 and 28.4]{H}.

\begin{proposition}
For each permutation $w$, the following hold:
\begin{enumerate}
\item The set $U_ww$ consists of the matrices \[\left\{w + u: \begin{array}{c} u \textup{ is nonzero only in entries that are } \\
\textup{{\bf both} above {\bf and} to the left of a nonzero entry in $w$} \end{array} \right\}.\]
\item The matrices $U_ww$ are a set of distinct coset representatives for the flags
in the Schubert cell $[Bw]$.  (See Figure \ref{uw}.)
\end{enumerate}
\end{proposition}

\begin{figure}[h]
\[  \left( \begin{array}{ccc}  a & 1 & 0 \\ 1 & 0 & 0 \\ 0 & 0 & 1 \end{array} \right)
\hspace{.5in} \left( \begin{array}{ccc}  a & b & 1 \\ 1 & 0 & 0 \\ 0 & 1 & 0 \end{array} \right)
\hspace{.5in}  \left( \begin{array}{ccc}  a & b & 1 \\ c & 1 & 0 \\ 1 & 0 & 0 \end{array} \right)\]
\caption{Examples of $U_ww$ when $n=3$ ($a$, $b$, and $c$ are free)} \label{uw}
\end{figure}

For each permutation $w$, let $w = s_{i_1,i_1+1}s_{i_2,i_2+1}\cdots s_{i_k,i_k+1}$ be a factorization with $k$ as small as possible.   We call $k$ the {\em length} of $w$, denoted
$\ell(w)$.  The length of $w$
relates the geometric, algebraic, and combinatorial descriptions of $GL_n/B$.  

\begin{proposition}
For each permutation $w$, the following hold:
\[ \begin{array}{rcl}
\ell(w) &=& \dim(\overline{[Bw]}) \\
&=& \textup{the number of nonzero
	entries (strictly) above the diagonal in }U_w\\
&=& \textup{the number of pairs } i<j \textup{ such that } w^{-1}(i)>w^{-1}(j).
\end{array}\]
%If adding Bott-Samelson resolution, add that length is also the minimum number of steps
%climbing up Bott-Samelson resolution chain
\end{proposition}

Each pair $i<j$ that satisfies $w^{-1}(i)>w^{-1}(j)$ is called an {\em inversion} for $w$.

The symmetric group is partially ordered by the {\em Bruhat order}.   The geometric
definition is that $v \leq w$ if and only if $[Bv] \subseteq \overline{[Bw]}$.
Combinatorially, we say $v \leq w$ if and only if there is a factorization 
$w = s_{i_1,i_1+1}s_{i_2,i_2+1}\cdots s_{i_k,i_k+1}$ so that
$v$ can be written as the product of a substring of the $s_{i_j,i_j+1}$.

\subsection{Properties of permutations and Schubert cells}

Several lemmas that follow from these properties will be used later
in this paper.  The difficulty of the proofs
depends on which characterization of the 
flag variety is used.  (Each of them is a nice exercise for the reader!)  

\begin{lemma} \label{combi lemma}
Fix $j  < k \leq n$.  For each permutation $w$, the following hold:
\begin{enumerate}
\item \label{shrinking w} The permutation $ws_{j,j+1}$ satisfies 
	$\ell(ws_{j,j+1}) = \ell(w)-1$ if
	and only if $w(j) > w(j+1)$.  Otherwise $\ell(ws_{j,j+1})=
		\ell(w)+1$.
\item \label{shrinking w by many}
	If $w(j) > w(k)$ then $w > ws_{jk}$ in the Bruhat order.
\end{enumerate}
\end{lemma}

\begin{proof}
The first part is classical, 
proven by noting that the sets of inversions of $w$ and of $ws_{j,j+1}$ differ
exactly by $(w(j), w(j+1))$.
%One proof compares the column of $w$ and of $ws_{j,j+1}$.
%This can also be proven by comparing the dimensions of $Uw$ and $Uws_{j,j+1}$.

To prove the next part, we 
show the closure of $[Bw]$ contains $ws_{jk}$.  Let 
$u_{w(k),w(j)}(a)$ be the upper-triangular matrix with $a$ in
position $(w(k),w(j))$, ones on the diagonal, and zeroes elsewhere.  Figure \ref{schematic of unipotent} is a schematic of $u_{w(k),w(j)}(a)w$.

\begin{figure}[h] \[ \left( \begin{array}{cc|c|ccc|c|cc}
	 \multicolumn{2}{c|}{\hspace{2em}} & \vdots & \multicolumn{3}{|c|}{} & \vdots & 
		\multicolumn{2}{|c}{\hspace{2em}} \\
	 \multicolumn{2}{c|}{\hspace{2em}} & 0 & \multicolumn{3}{|c|}{} & 0 & 
		\multicolumn{2}{|c}{\hspace{2em}} \\
	\cline{1-9} \multicolumn{2}{c|}{0 \cdots 0} & a & \multicolumn{3}{|c|}{0 \cdots 0} & 1 & 
		\multicolumn{2}{|c}{0 \cdots 0} \\
	\cline{1-9} \multicolumn{2}{c|}{} & 0 & \multicolumn{3}{|c|}{} & 0 & \multicolumn{2}{|c}{} \\
	 \multicolumn{2}{c|}{} & \vdots & \multicolumn{3}{|c|}{} & \vdots & \multicolumn{2}{|c}{} \\
	 \multicolumn{2}{c|}{} & 0 & \multicolumn{3}{|c|}{} & 0 & \multicolumn{2}{|c}{} \\
	\cline{1-9} \multicolumn{2}{c|}{0 \cdots 0} & 1 & \multicolumn{3}{|c|}{0 \cdots 0} & 0 & 
	\multicolumn{2}{|c}{0 \cdots 0} \\
	\cline{1-9} \multicolumn{2}{c|}{} & 0 & \multicolumn{3}{|c|}{} & 0 & \multicolumn{2}{|c}{} \\
	 \multicolumn{2}{c|}{} & \vdots & \multicolumn{3}{|c|}{} & \vdots & \multicolumn{2}{|c}{}  \\
	 \multicolumn{2}{c|}{} & 0 & \multicolumn{3}{|c|}{} & 0 & \multicolumn{2}{|c}{} 
	 	\end{array}\right) \]
		\caption{Schematic of $u_{w(k),w(j)}(a)w$} \label{schematic of unipotent}
		\end{figure}

Denote the flag
$[u_{w(k),w(j)}(a)w]$ by $V_1 \subseteq
V_2 \subseteq \cdots \subseteq V_n$.  
The first $j-1$ subspaces and last $n-k+1$ subspaces of 
this flag agree with those of the flag $[ws_{jk}]$.  
The other subspaces are
\[\begin{array}{l}\langle V_{j-1}, e_{w(j)}+ae_{w(k)} 
	\rangle \subseteq \langle V_{j-1}, e_{w(j)}+ae_{w(k)},
e_{w(j+1)} \rangle \subseteq \cdots \\
\hspace{.5in} \subseteq \langle V_{j-1}, e_{w(j)}+ae_{w(k)} ,
e_{w(j+1)}, e_{w(j+2)}, \ldots, e_{w(k-1)} \rangle  \\
\hspace{.5in} \subseteq
\langle V_{j-1}, e_{w(j)},e_{w(j+1)}, e_{w(j+2)}, \ldots, e_{w(k-1)}, e_{w(k)} \rangle = V_{k}.
\end{array}\]    
As $a$ approaches $\infty$, these subspaces approach the subspaces
\[\begin{array}{l}\langle V_{j-1}, e_{w(k)} \rangle \subseteq \langle V_{j-1}, e_{w(k)} ,
e_{w(j+1)} \rangle \subseteq \cdots \\
\hspace{.5in} \subseteq \langle V_{j-1}, e_{w(k)},
e_{w(j+1)}, e_{w(j+2)}, \ldots, e_{w(k-1)} \rangle   \\
\hspace{.5in} \subseteq
\langle V_{j-1}, e_{w(k)},e_{w(j+1)}, e_{w(j+2)}, \ldots, e_{w(k-1)}, e_{w(j)} \rangle = V_k,
\end{array}\]    
which are the corresponding
parts of $[ws_{jk}]$.  Thus 
$\lim_{a \mapsto \infty} [u_{w(k), w(j)}(a)w]=[ws_{jk}]$.
It follows that $[Bws_{jk}] \subseteq \overline{[Bw]}$, and so
$ws_{jk} < w$.
\end{proof}

\subsection{Hessenberg varieties} \label{hess vars}

In this section, we define Hessenberg varieties algebraically.  
We also discuss some technical issues that arise.

To obtain an algebraic characterization of Hessenberg varieties, we use subspaces of
$n \times n$ matrices rather than the Hessenberg function $h$.  The matrix basis
unit that is zero except in entry $(i,j)$, where it is one, is denoted $E_{ij}$.
Each Hessenberg function defines a subspace of $n \times n$
matrices by $H_h = \langle E_{ij}: i \leq h(j) \rangle$.
We call $H_h$ a Hessenberg space.  
The Hessenberg variety of $X$ and $h$ is
\[\He(X,h) = \{ \textup{Flags }[g]: g^{-1}Xg \in H_h\}.\]

Many examples of Hessenberg spaces come from classical Lie theory.
If $h$ is the Hessenberg function with $h(i)=i$ for each $i$ then $H_h$ is the set of 
upper-triangular matrices.  If $h$ is the Hessenberg function given
by $h(i)=n$ for each $i$ then $H_h$ consists of all $n \times
n$ matrices.  In fact, if $H_h$ is any parabolic subalgebra, then $H_h$ is a Hessenberg space and the corresponding $h$ is one of the parabolic Hessenberg functions from the Introduction.

Most Hessenberg spaces are not parabolic.  For instance, the Hessenberg function given
by $h(i)=i+1$ when $i \neq n$ and $h(n)=n$ corresponds to the subspace
$H_h$ which is zero below the subdiagonal.  Figure
\ref{Hess 1 diagram} shows this for $n=4$.  
\begin{figure}[h]
\[\begin{array}{rcl} h(1) & = & 2 \\ h(2) & = & 3 \\ h(3) & = & 4 \\ h(4) & = & 4 \end{array}
\longleftrightarrow
\left(\begin{array}{cccc} * & * & * & *  \\  {*} & * & * & * \\ 
0 & * & * & * \\ 0 & 0 & * & * \end{array} \right)\]
\caption{One Hessenberg function and space when $n=4$} \label{Hess 1 diagram}
\end{figure}
Hessenberg varieties with this Hessenberg function are important in various
applications, including numerical analysis \cite{dMPS} and computing quantum 
cohomology of the flag variety (see \cite{K} and \cite{R}).

Our definition of Hessenberg functions 
omits one condition from the original definition in \cite{dMPS}, which
 also requires $h(i) \geq i$ for each $i$.  
This paper studies a strictly larger collection of varieties than in \cite{dMPS}.
Our generalization is particularly useful when 
$X$ is nilpotent.  (When $X$ is regular semisimple, the variety $\He(X,h)$ will be
empty if $h(i)<i$ for each $i$.)  Nilpotent Hessenberg varieties
arise naturally when studying representations of the 
symmetric group on Hessenberg varieties that generalize Springer's
correspondance \cite{MT}.

Section \ref{generalizations}
generalizes this definition (and other results) to all Lie types.

Our first proposition establishes that nilpotent Hessenberg varieties 
depend only on the $i$ for which the Hessenberg function 
does not satisfy $h(i)=i$.  

\begin{proposition} \label{diagonals are equivalent}
Fix $n$ and fix $i$ such that $1 \leq i \leq n$.  Suppose $h$ is a Hessenberg
function with $h(i)=i$ and that the function $h'$ defined by 
\[h'(j)= \left \{ \begin{array}{rl} h(j) & \textup{ if }j \neq i, \textup{ and} \\ 
i-1 &\textup{ for } i=j \end{array} \right.\] is also a Hessenberg function.
If $X$ is nilpotent then $\He(X,h)=\He(X,h')$.
\end{proposition}

\begin{proof}
If $g^{-1}Xg \in H_{h'}$ then $g^{-1}Xg \in H_h$ since $H_{h'} \subseteq H_h$.  
Now assume $g^{-1}Xg\in H_h$. 
We have
$$
(g^{-1}Xg)  e_i \in c e_i +  \langle e_1,\ldots, e_{i-1} \rangle
$$  
where $e_j$ are the standard basis vectors
for $\C^n$. Also 
$$
(g^{-1}Xg) \langle e_1,\ldots, e_{i-1}\rangle \  
\subseteq\  \langle e_{1},\ldots, e_{h(i-1)}\rangle  \ 
\subseteq \ \langle e_1,\ldots, e_{i-1}\rangle, 
$$
since $h(i-1)<h(i)=i$. Since $X$ is nilpotent, applying $g^{-1}Xg$ to $e_i$ sufficiently many 
(e.g. $n$) times should
give zero. On the other hand we have
$$
(g^{-1}X{g})^n e_i \in c^n e_i +  \langle e_1,\ldots, e_{i-1}\rangle. 
$$
Therefore $c=0$, and as a consequence $g^{-1}Xg$ lies in $H_{h'}$. 
%For each upper-triangular $b$, the conjugate $b^{-1}Xb$ can be written
%$b^{-1}Xb=\sum c_{ij}E_{ij}$, 
%since $X$ is strictly upper-triangular.  For each permutation
%$w$, the conjugate $w^{-1}(\sum c_{ij}E_{ij})w = \sum c_{ij}E_{w^{-1}(i),w^{-1}(j)}$.
%In particular, $w^{-1}b^{-1}Xbw$ is zero along the diagonal for each 
%upper-triangular matrix $b$ and permutation $w$.  Since $h$ and $h'$ differ only
%by entries along the diagonal, the Hessenberg varieties $\He(X,h) = \He(X,h')$.
\end{proof}

Comments from K.\ Rietsch greatly improved this proof.
This lemma motivates the following definition, also suggested by K.\ Rietsch.

\begin{definition}
For each linear operator $X$, the Hessenberg spaces $H$ and $H'$ are $X$-equivalent if
$\He(X,H) = \He(X,H')$.   In this case, we write $H \sim_X H'$ and say that $H$ and $H'$ are in the same $X$-equivalence class.
\end{definition}

$X$-equivalence of Hessenberg functions is defined the same way.

The $X$-equivalence class of Hessenberg spaces (or functions) depends only on the 
conjugacy class of $X$, since $\He(X,H) \cong \He(g^{-1}Xg,H)$ (see \cite[Proposition 2.7]{T}).

For instance, if $X=0$ then there is only one $X$-equivalence class of Hessenberg spaces.  
If $X$ is nilpotent, then the Hessenberg function defined by $h(i)=i$ for all $i$ is $X$-equivalent to
the function defined by $h'(i)=i-1$ for all $i$.  Alternatively,
 the Hessenberg space consisting of all upper-triangular matrices is $X$-equivalent to the space
of all strictly upper-triangular matrices.  (This fact is used frequently in Springer theory.)  We generalize this in the next corollary, whose proof is immediate from Proposition \ref{diagonals are equivalent}.

\begin{corollary}
For each nilpotent linear operator $X$, there is a unique minimal element of each 
$X$-equivalence class
of Hessenberg functions (respectively Hessenberg spaces).  This minimal element satisfies
\begin{itemize}
\item if there exists $i$ such that $h(i)=i$, then $h(i-1)=i$ as well;
\item if there exists a matrix $\sum c_{jk}E_{jk}$ in $H_h$ and $i$ such that
	the coefficient $c_{ii} \neq 0$, then $E_{ii}$ and $E_{i-1,i}$ are both in $H_h$.
\end{itemize}
\end{corollary}

Typically, we assume $H$ and $h$ are minimal in their $X$-equivalence classes.

\section{Geometry and topology of $\X{h}$} \label{highweight}

In this section,
we fix $X$ to be the matrix $E_{1n}$ and study the Hessenberg varieties
\[\X{h} = \{ \textup{Flags } [g]: g^{-1} E_{1n} g \in H_h \} = \{ \textup{Flags } V_1 \subseteq \cdots 
     \subseteq V_n: E_{1n} V_i  \subseteq V_{h(i)} \}.\]
We will show that these Hessenberg varieties are unions of Schubert varieties.
Loosely speaking, each Schubert variety comes from one ``corner" of the 
Hessenberg space.  We will identify explicitly these Hessenberg varieties,
including which Schubert varieties arise and their dimensions.  We will
also show that many of these Hessenberg varieties are not pure-dimensional.

% proof in general is $u^{-1}X_{\theta}u = X_{\theta}$
\begin{proposition} \label{schubert varieties}
$\X{h}$ is a union of Schubert varieties $\bigcup \Cy_w$.
\end{proposition}

\begin{proof}
Each flag can be written in row echelon form as $[uw]$ 
for some invertible upper-triangular $u$ and permutation matrix $w$.  
The flag $[uw]$ is in $\X{h}$ if and only if 
$w^{-1}u^{-1}E_{1n}uw$ is in $H_h$.  Direct calculation shows that
$u^{-1}E_{1n}u$ is a nonzero scalar multiple of $E_{1n}$ for each upper-triangular
$u$.  Thus, the flag $[uw]$ is in $\X{h}$ if and only if $[w]$ is in $\X{h}$.

This means $\X{h}$ is a union of Schubert cells, say $\X{h} = \bigcup [Bw]$,
and so $\X{h} \subseteq \bigcup \Cy_w$.  
Since $\X{h}$ is closed, it also
contains the closures $\bigcup \overline{[Bw]} =  \bigcup \Cy_w$.
\end{proof}

In general the variety $\He(X,h)$ is not a union of Schubert cells $[Bw]$.
In fact, if $g^{-1}Xg$ is another element of the conjugacy class of $X$, then
typically at most one of $\He(X,h)$ and $\He(g^{-1}Xg,h)$ is a union of 
cells $[Bw]$, even though
the two varieties
are homeomorphic \cite[Proposition 2.7]{T}.  For instance,  
suppose $n=3$ and the Hessenberg function satisfies
$h(i)=i$ for each $i$.  
Each of $\He(E_{12},h)$ and $\He(E_{13},h)$ is homeomorphic to two copies of
$\Pp^1$ glued together at a point.
%, as in Figure \ref{subregular Springer}.  
However, the variety 
$\He(E_{13},h)$ is the union of the Schubert varieties $\Cy_{s_1} \cup
\Cy_{s_2}$, while $\He(E_{12},h)$ is a one-dimensional
closed subvariety of $\Cy_{s_2 s_1}$.

\begin{figure}[h]
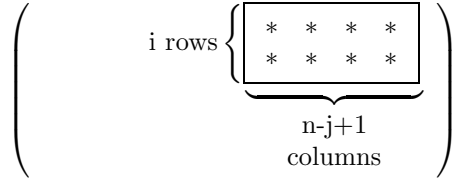

\[\left( \begin{array}{ccr}
\hspace{1em} & \hspace{1em} & \mbox{i rows} \left\{\fbox{$\begin{array}{cccc} * & * & * & * 
	\\ {*} & * & * & {*} \end{array}$}\right.  \vspace{-.7em} \\ 
	 & & \underbrace{\makebox[6.7em]{}}_{\mbox{\begin{tabular}{c}n-j+1 \\ columns \end{tabular}}}
\end{array} \right)\]
\caption{Schematic diagram of $H_{ij}$} \label{Hij}
\end{figure}

%\begin{figure}[h]
%\begin{todo}
%add picture!
%\end{todo}
%\caption{The subregular Springer fiber when $n=3$} \label{subregular Springer}
%\end{figure}
For each $i \neq j$, let $h_{ij}$ be the Hessenberg function defined by 
\[h_{ij}(k) = \left\{ \begin{array}{ll} 0 & \textup{ if } k < j \textup{ and} \\
        i & \textup{ if } k \geq j, \end{array} \right.\]
The corresponding Hessenberg space $H_{ij}$ is spanned by the matrix basis units
$E_{kl}$ with $k \leq i$ and $l \geq j$.  In other words, $H_{ij}$ is the subspace of matrices
which are zero outside of the upper-right $i \times (n-j+1)$ rectangle, as in Figure
\ref{Hij}.
For example, $H_{n1}$ consists of all $n \times n$ matrices and $H_{1n}$ is
just the span of $E_{1n}$.  If the sun rises at the far left of the $i^{th}$ row,
travels around the bottom left corner of the matrix, and 
sets at the bottom of the $j^{th}$ column, then $H_{ij}$ is the shadow cast by the matrix
basis unit $E_{ij}$ during the course of the `day'.  (A.\ Ottazzi created this image in \cite{O}.)

\begin{lemma} \label{main schubert lemma}
For each pair $i \neq j$, let $w$ be the permutation that has 
$e_n$ in column $j$, $e_1$ in column $i$, and the other vectors inserted
in decreasing order $(e_{n-1}$, $e_{n-2}$, $\ldots$, $e_2)$ in the remaining columns.
Then $\X{H_{ij}} = \Cy_w$.
\end{lemma}

\begin{proof}
The proof has three parts.  First, we show that if $s$ is a permutation, then 
$[s] \in X_{H_{ij}}$ if and only if 
$s$ has $e_1$ somewhere in its first $i$ columns and $e_n$ 
 in its last $n-j+1$ columns.  For each such $s$, we form the permutation $s'$ 
 by moving $e_1$ to the $i^{th}$
column, moving $e_n$ to the $j^{th}$ column, and keeping the other columns in
the same order as in $s$.  We then show that $s' \geq s$.  Finally, we show that $w \geq s'$.

Suppose the matrix $s$ has $e_1$ in its $k^{th}$ column and $e_n$ in its $l^{th}$ column.  
Since
$s^{-1} = s^t$, the matrix $s^{-1}$ has $e_1$ in its $k^{th}$ row.  So  
\begin{equation} \label{max root} s^{-1}E_{1n}s = E_{kl}. \end{equation}  
This is in $H_{ij}$ if and only if $k \leq i$ and $l \geq j$.  
We conclude that the flag $[s] \in X_{H_{ij}}$ if and only if
$s^{-1}e_1 = e_k$ for  $k \leq i$ and 
$s^{-1}e_n = e_l$ for $l \geq j$.

The permutation $w$ satisfies
this condition so $\X{H_{ij}} \supseteq \Cy_w$.  
We now show that for any permutation $s$ of this form, the flag 
$[s]$ is in the variety $\Cy_w$.

We begin by moving the column with $e_n$ to the left or  
the column with $e_1$ to the right, as long as one of those
moves is possible.  Suppose $l > j$ and either the $(l-1)^{th}$ column is not
$e_1$ or it is $e_1$ and $l-1 \neq i$.  The flag $[ss_{l-1,l}]$
is also in $\X{H_{ij}}$.  Lemma \ref{combi lemma} Part
\ref{shrinking w by many} showed 
$ss_{l-1,l} > s$, so the corresponding Schubert varieties satisfy 
$\Cy_s \subseteq \Cy_{ss_{l-1,l}}$.  
(When $k < i$ and either the $(k+1)^{th}$ column is not $e_n$ or
it is $e_n$ but $k+1 \neq j$, use the flag $[ss_{k,k+1}]$ in a symmetric argument.)  

A move of this sort will be impossible exactly when $j<i$ and either 
\vspace{.05in} \begin{itemize}
\item $i=k$ and $l=k+1=i+1$ or
\item $l=j$ and $k=j-1=l-1$. \vspace{-.5in}
\end{itemize}

\vspace{.1in} \hspace{2.5in}  $\left( \begin{array}{rrrr} \hspace{1em} \cdots & \begin{tabular}{|c|}\cline{1-1} \\
$e_1$ \\ \\ \hline \end{tabular} & \begin{tabular}{|c|}\cline{1-1} \\
$e_n$ \\ \\ \hline \end{tabular} & \cdots \hspace{1em} \end{array} \right)$

The diagram is a schematic for these cases: the vectors $e_1$ and $e_n$
are adjacent, and the $i^{th}$ column is in place (respectively $j^{th}$) while $e_n$
is moving to the left (respectively $e_1$ to the right).
Lemma \ref{combi lemma} Part \ref{shrinking w by many} 
shows that the permutation
obtained from $s$ by exchanging its $(i+1)^{th}$ and $(i-1)^{th}$ columns
is greater than $s$ in the Bruhat order (respectively $j-1$ and $j+1$).

Once $e_1$ is to the right of $e_n$, successively multiply $s$ on the right 
by $s_{k,k+1}$ or $s_{l-1,l}$ to obtain
a permutation $s'$ with $s' \geq s$, so that $s'(e_i) = e_1$ and $s'(e_j) = e_n$.
%$i^{th}$ column of $s'$ is $e_1$ and $j^{th}$ column of $s'$ is $e_n$

We now prove by induction that $s' \leq w$.  Assume that the
first $t$ columns of $s'$ and $w$ agree, and the $(t+1)^{th}$ does not.  
The $(t+1)^{th}$ column of $w$ is
filled with $e_{w(t+1)}$. 
Neither $s'(t+1)$ nor $w(t+1)$ is
in $\{1,n\}$ because $s'(t+1) \neq w(t+1)$.
Since $w$ and $s'$ agree in the first $t$ columns, the
column vector $w(e_{t+1})$ is none of $e_{s'(1)}$, $e_{s'(2)}$, $\ldots$, $e_{s'(t)}$,
so 
%there is a positive integer $t_1$ with $s'(e_{t+1})=e_{w(t+1+t_1)}$.  
there is a positive integer $t_1$ such that
$s'(e_{t+1+t_1})=e_{w(t+1)}$.   The permutation $s'' = s' s_{t+1, t+1+t_1}$ 
satisfies $s'' \geq s'$ by Lemma \ref{combi lemma} Part 
\ref{shrinking w by many}.  
Since neither $e_1$ nor $e_n$ moved, $s''$ has
$s''(e_i)=e_1$ and $s''(e_j)=e_n$, and also 
agrees with $w$ in its first $t+1$ columns.  By 
induction, the claim follows.
\end{proof}

The following corollary restates the condition on $w$.

\begin{corollary}
For each pair $i \neq j$, let $w$ be the largest permutation in the Bruhat order that
satisfies $w^{-1}E_{1n}w = E_{ij}$.  Then $X_{H_{ij}} = \Cy_w$.
\end{corollary}

We can factor $w$ explicitly in terms of simple transpositions.

\begin{corollary} \label{factoring w}
Let $w_0$ be the permutation with $w_0e_k=e_{n-k+1}$ for each $k=1, \ldots, n$.
For each pair $i \neq j$, 
the Hessenberg variety $X_{H_{ij}}= \Cy_w$, where
\[w = \left\{ \begin{array}{ll} 
w_0s_{12} s_{23} \cdots s_{j-1,j}s_{n,n-1} \cdots s_{i+1,i} & \textup{ if } j<i \textup{ and} \\
w_0s_{12} s_{23} \cdots s_{j-2,j-1}s_{n,n-1} \cdots s_{i+1,i} & \textup{ if } j>i.
\end{array} \right.\]
\end{corollary}

\begin{proof}
For each matrix $M$,
the product $Ms_{12}s_{23}\cdots s_{k,k+1}$ cyclically permutes the first $k+1$ columns of 
$M$, sending the first column to the $(k+1)^{th}$ position and moving each of the
other columns one position to the left.  Similarly, the product $Ms_{n,n-1}s_{n-1,n-2} \cdots 
s_{k+1,k}$ cyclically permutes the last $n-k+1$ columns, moving the last column to the
$k^{th}$ and moving the others one column to the right.  Cyclically permuting the first
$j$ (respectively $j-1$) 
columns and the last $n-i+1$ columns of $w_0$ gives the permutation $w$ of
Lemma \ref{main schubert lemma}.
\end{proof}

This gives a closed formula for the dimension of $X_{H_{ij}}$.

\begin{corollary} \label{closed formula for dimn}
For each $i \neq j$, the dimension of $X_{H_{ij}}$ is 
\[\left\{ \begin{array}{ll}
\vspace{.5em} \left( \begin{array}{c} n \\ 2 \end{array} \right) - \left(j-1+n-i
\right) & \textup{if }j<i \textup{ and } \\
\left( \begin{array}{c} n \\ 2 \end{array} \right) - \left(j-2+n-i
\right) & \textup{if }j>i. \end{array} \right.\]
\end{corollary}

\begin{proof}
The length of the permutation $w_0$ is $\left( \begin{array}{c} n \\ 2 \end{array} \right)$.
Let $w = w_0 \prod s_{k,k+1}$ be the factorization from Corollary \ref{factoring w}.
Each simple transposition in this factorization reduces the length of $w_0$ by one, from
Lemma \ref{combi lemma} Part \ref{shrinking w}.
\end{proof}

\subsection{The components of $\X{H}$}

It is usually difficult to identify the irreducible components of Hessenberg varieties.
However, when $X = E_{1n}$, it can be done.

\begin{proposition} \label{decomp into unions}
For all $H$ and $H'$, we have $\X{H \cup H'} = \X{H} \cup \X{H'}$.
\end{proposition}

\begin{proof}
The flag $[w]$ is in $\X{H \cup H'}$ if and only if $w^{-1}E_{1n}w$ is in $H \cup H'$.  Since
$w^{-1}E_{1n}w$ is a matrix basis unit, it is in $H \cup H'$ if and only if
either $w^{-1}E_{1n}w$ is in $H$ or $w^{-1}E_{1n}w$ is in $H'$.  This holds
if and only if the flag $[w]$ is in $\X{H} \cup \X{H'}$.
\end{proof}

\begin{lemma} \label{containment lemma}
Let $H$ and $H'$ be Hessenberg spaces that are minimal in their $E_{1n}$-equivalence classes.
Then $\X{H} \subseteq \X{H'}$ if and only if $H  \subseteq H' $.
\end{lemma}

\begin{proof}
We reduce to the case when $\X{H}$ and $\X{H'}$ are Schubert varieties.  Write
$H = \bigcup H_{ij}$ and $H' = \bigcup H_{i'j'}'$.    Each $H_{ij}$ satisfies $i \neq j$
since $H$ is minimal 
in its $E_{1n}$-equivalence class (respectively $i' \neq j'$).
For each pair $i \neq j$, the matrix
$E_{ij}$ is in $H'$ if and only if $E_{ij} \in H_{i'j'}'$ for some $i',j'$.  This holds 
if and only if $H_{ij} \subseteq H_{i'j'}'$.  Consequently 
$H \subseteq H'$ if and only if for each $i,j$ there exist $i',j'$ such that 
$H_{ij} \subseteq H_{i'j'}'$.  We know
$\X{H} = \bigcup \X{H_{ij}}$ and $\X{H'} = \bigcup \X{H_{i'j'}'}$ from 
Proposition \ref{decomp into unions}.  It suffices to show that 
$\X{H_{ij}} \subseteq \X{H_{i'j'}'}$ if and only if $H_{ij} \subseteq H_{i'j'} $.
 
Both $\X{H_{ij}}$ and $\X{H_{i'j'}}$ are a disjoint union of Schubert cells
by Proposition \ref{schubert varieties}.  This means 
the inclusion $\X{H_{ij}} \subseteq \X{H_{i'j'}}$ holds if and only if each permutation
flag $[s]$ in $\X{H_{ij}}$ is also contained in $\X{H_{i'j'}}$.  Equation 
\ref{max root} shows that $[s]$ is in $\X{H_{ij}}$ if and only if $s^{-1}E_{1n}s = E_{kl}$, where
$k$ and $l$ satisfy the conditions $k \leq i$ and $l \geq j$.   It follows that
each permutation flag $[s]$ in $\X{H_{ij}}$ is also in $\X{H_{i'j'}}$ if and only if
$i \leq i'$ and $j \geq j'$, which is true if and only if $H_{ij} \subseteq H_{i'j'}$.
\end{proof}

\begin{definition}
A maximal decomposition of the Hessenberg space 
$H$ is a union $H = \bigcup H_{ij}$ so that no pair $H_{ij}$, $H_{i'j'}$
satisfies $H_{ij} \subseteq H_{i'j'}$.
\end{definition}

If $H$ is minimal in its $E_{1n}$-equivalence class, then a maximal decomposition
$H = \bigcup H_{ij}$ further satisfies $i \neq j$ for each $H_{ij}$.

\begin{corollary} \label{identify components}
Let $H$ be minimal in its $E_{1n}$-equivalence class.
If $H = \bigcup H_{ij}$ is a maximal decomposition,
 the components of $\X{H}$ are the Schubert varieties $\X{H_{ij}}$.
\end{corollary}

\begin{proof}
Write $\X{H} = \bigcup \X{H_{ij}}$ as in Proposition \ref{decomp into unions}.
For each $H_{ij}$, there is a unique permutation $w_{ij}$ such that $[Bw_{ij}]$ is dense in
$\X{H_{ij}}$ by Lemma \ref{main schubert lemma}.  For every $(i',j') \neq (i,j)$, 
Lemma \ref{containment lemma}
shows that $[w_{i'j'}]$ is not in 
$\X{H_{ij}}$, and so $[Bw_{i'j'}] \cap \X{H_{ij}}$ is empty.
This means $\X{H_{ij}}$ is an irreducible component of $\X{H}$.
\end{proof}

\begin{corollary} \label{banded form}
Fix $H$, a minimal Hessenberg space in its $E_{1n}$-equivalence class.
The Hessenberg variety $\X{H}$ is pure dimensional if and only if there exists 
an integer $k \in 
\{1, 2, \ldots, n-1\}$ and a subset $I \subseteq \{1, 2, \ldots, n-k\}$ such that
either $H = \bigcup_{i \in I} H_{i,i+k}$ or $H= \bigcup_{i \in I} H_{i,i-k}$.
\end{corollary}

\begin{proof}
Let $H = \bigcup H_{ij}$ be a maximal decomposition and 
write $\X{H}$ as a union of its irreducible components $X_{H_{ij}}$.
Each $X_{H_{ij}}$ is a Schubert variety that has
dimension $\left( \begin{array}{c} n \\ 2 \end{array} \right) - \left(j-1+n-i
\right)$ if $j<i$ and $\left( \begin{array}{c} n \\ 2 \end{array} \right) - \left(j-2+n-i
\right)$ if $j>i$ by Corollary \ref{closed formula for dimn}.  Given pairs $(i,j)$ and $(i',j')$,
the varieties $X_{H_{ij}}$ and $X_{H_{i'j'}}$ have the same dimension
if and only if $j-i = j'-i'$.
\end{proof}

This gives a collection of examples of nilpotent Hessenberg varieties that are connected and 
not pure dimensional.   For instance, the Hessenberg space
$H = H_{41} \cup H_{54}$ of
$5 \times 5$ matrices gives a variety $\X{H}$ in $GL_5/B$ that is not pure.

\section{A semisimple Hessenberg variety that is not pure dimensional} \label{semisimple impurity}

In this section, we describe another way that Hessenberg varieties can fail to be pure 
dimensional.  
The next proposition generalizes an example that R.\ MacPherson and I discovered.

\begin{proposition}
Fix $X = \sum_{i=1}^{n-1} E_{ii}$.  Let $h$ be the Hessenberg function with 
$h(i) = n-1$ for all $i \leq n-1$, and $h(n)=n$.
The variety $\He(X,h)$ is the disjoint union of two components,
one of which is homeomorphic to $GL_{n-1}/B$ and the other of which is
homeomorphic to a fiber bundle over $\Pp^{n-2}$ with fiber $GL_{n-1}/B$.  
In particular, the Hessenberg variety $\He(X,h)$ is not pure dimensional.
\end{proposition}

\begin{proof}
By definition, each flag $V_{\bullet}$ in $\He(X,h)$ satisfies $X V_{n-1} \subseteq
V_{n-1}$.  Since $X (\sum_{i=1}^n a_i e_i) = \sum_{i=1}^{n-1} a_i e_i$, 
either
\begin{enumerate}
\item \label{first component} $e_n \in V_{n-1}$ or
\item \label{second component} $V_{n-1} = \langle e_1, e_2, \ldots, e_{n-1} \rangle$.
\end{enumerate}
These conditions are closed and so define two closed subvarieties ${\mathcal Y}_1$ 
and ${\mathcal Y}_2$, respectively,
in $GL_n/B$.  The two conditions cannot be simultaneously satisfied so
$\He(X,h)$ is the disjoint union ${\mathcal Y}_1 \cup {\mathcal Y}_2$.  We now describe these
subvarieties.

First we show that ${\mathcal Y}_2 \cong GL_{n-1}/B$.  The flag $[g]$ satisfies
Condition \ref{second component} if and only if the matrix $g$ is in 
\[P = \left( \begin{array}{cccc|c} \multicolumn{4}{c|}{ GL_{n-1}} & \begin{array}{c} {*} \\ {*}
	\\ \vdots \\ {*} \end{array} \\
	 \cline{1-5} 0 & 0 & \cdots & 0 & \C^* \end{array} \right).\]  
	 In other words,
the component ${\mathcal Y}_2$ is isomorphic to $GL_{n-1}/B$ via the isomorphism
that sends $V_1 \subseteq \cdots \subseteq V_{n-1} \subseteq V_n$ to the flag
$V_1 \subseteq \cdots \subseteq V_{n-1}$ inside $\langle e_1, \ldots, e_{n-1} \rangle$.

Now we study ${\mathcal Y}_1$.    Denote the Grassmannian of $n-1$-planes in $\C^n$ by
$G(n-1,n)$.  Write $\pi_{n-1}: GL_{n}/B \longrightarrow G(n-1,n)$ for the projection
that sends the flag $V_1 \subseteq V_2 \subseteq \cdots \subseteq V_n$ to the subspace $V_{n-1}$.
This is a continuous map; in fact, it is the quotient map $\pi_{n-1}: GL_{n}/B \longrightarrow GL_{n}/P$.

Restrict the map to $\pi_{n-1}|_{{\mathcal Y}_1}: {\mathcal Y}_1 \longrightarrow G(n-1,n)$.  The image 
$\pi_{n-1}({\mathcal Y}_1)$ is 
\[ \pi_{n-1}({\mathcal Y}_1) = \{\textup{Subspaces } V_{n-1} \textup{ such that } 
	e_n \in V_{n-1}\}.\]
This is isomorphic to the set of $n-2$-dimensional subspaces in $\langle e_1, \ldots, e_{n-1} 
\rangle$, so $\pi_{n-1}({\mathcal Y}_1) \cong G(n-2,n-1)$.  Since $G(n-2,n-1) \cong
\Pp^{n-2}$, we conclude that the image $\pi_{n-1}({\mathcal Y}_1) \cong \Pp^{n-2}$.

We now identify the fiber $\left(\pi_{n-1}|_{{\mathcal Y}_1}\right)^{-1}(V_{n-1})$ of each $V_{n-1} \in 
\pi_{n-1}({\mathcal Y}_1)$.
The flag $W_1 \subseteq \cdots \subseteq W_n$ is in 
$\left(\pi_{n-1}|_{{\mathcal Y}_1}\right)^{-1}(V_{n-1})$ if and only
if $W_{n-1} = V_{n-1}$.  Every flag in $GL_n/B$ satisfies $W_n = \C^n$, so the fiber
is characterized by 
\[\left(\pi_{n-1}|_{{\mathcal Y}_1}\right)^{-1}(V_{n-1}) =  \{\textup{Flags such that } 
W_1 \subseteq W_2 \subseteq \cdots W_{n-2} \subseteq V_{n-1}\}.\]  
This is the
set of complete flags in $V_{n-1}$ and is homeomorphic to $GL_{n-1}/B$.

Consequently, the map $\pi_{n-1}: {\mathcal Y}_1 \longrightarrow \pi_{n-1}({\mathcal Y}_1)$ is
a fiber bundle whose base space is homeomorphic to $\Pp^{n-2}$ and whose fiber 
is homeomorphic to $GL_{n-1}/B$. 
\end{proof}

For example, when $n=3$ the Hessenberg variety
$\He(X,h)$ is a disjoint union of $\Pp^1$ and a $\Pp^1$-bundle over $\Pp^1$.

\section{Generalizing to all Lie types} \label{generalizations}

In this section, we discuss generalizations of these results to arbitrary Lie type.  Our 
exposition is brief; we assume our reader is familiar with the general theory.

Let $G$ be a complex reductive 
linear algebraic group, $\gp$ its Lie algebra, $B$ a fixed Borel subgroup,
and $\bo$ its Lie algebra.  The full flag variety is $G/B$ and its elements are
written $[g]$.  Let $T$ be a maximal torus contained in $B$ and $\tor$ be the
Cartan subalgebra associated to $T$.  We will also use $\mathfrak{n}^-$, the 
maximal nilpotent subalgebra in the opposite Borel subalgebra $\bo^-$.  
Let $W$ be the Weyl group.

The positive roots in the root system corresponding to $\gp$ are denoted $\Phi^+$ and the negative roots are $\Phi^-$.  The inner product on $\Phi$ is written $\langle \cdot, \cdot \rangle$.
We refer to the length of roots, which can be either short or long.
If $\alpha$ and $\beta$ are two roots, then $\alpha \succ \beta$ means $\alpha - \beta$ is a 
{\em sum}
of positive roots.  (Note that this is not the partial ordering where $\alpha > \beta$ means $\alpha-
\beta$ is a positive root.)  If $\alpha=\sum c_i \alpha_i$ is a (reduced) sum of simple roots,
then the support of $\alpha$ is the set $\textup{supp}(\alpha) = \{ \alpha_i: c_i \neq 0 \}$.
Given $\alpha$, we write $E_{\alpha}$ for a root vector corresponding to $\alpha$.

A Hessenberg space $H$ is a linear subspace of  matrices such
that $[H,\bo] \subseteq H$.  (This definition omits one condition from that found in 
\cite{dMPS}.)  Suppose $X$ is in $\gp$ and $H$ is a Hessenberg 
space.  The Hessenberg variety of $(X,H)$ is given by
\[\He(X,H) = \{ [g] \in G/B : g^{-1}Xg \in H\}.\]

\begin{proposition}
Let $E_{\theta}$ be a weight
vector for the highest weight $\theta$.  For each Hessenberg space $H$, the variety $\He(E_{\theta},H)$ is a union of  Schubert varieties.
\end{proposition}

This generalizes Proposition \ref{schubert varieties}.  The proof is the same as in 
Proposition \ref{schubert varieties}: the flag $[bw]$ is in $\He(E_{\theta}, H)$ if and only if
$w^{-1}b^{-1}E_{\theta}bw$ is in $H$, and 
 the adjoint action of $B$ multiplies $E_{\theta}$ by a nonzero constant factor.  

\begin{definition}
For each root $\alpha$, define $H_{\alpha}$ to be minimal with respect to inclusion
among all Hessenberg
spaces that contain the root vector $E_{\alpha}$.
\end{definition}

If $\alpha$ is positive then $H_{\alpha}$ is 
the span of the root vectors $E_{\beta}$ with $\beta \succeq \alpha$.
However, this is not true when $\alpha$ is negative.  In that case, 
every positive root $\beta$ satisfies $\beta \succeq \alpha$,
but $H_{\alpha}$ need not contain $\bo$.  

Let $N(\textup{supp}(\alpha))= \{\alpha_j: \exists \alpha_i \in \textup{supp}(\alpha) \textup{ with } 
\langle \alpha_j, \alpha_i \rangle \neq 0\}$.  In other words, $N(\textup{supp}(\alpha))$ consists
of $\textup{supp}(\alpha)$ as well as the simple roots that are
joined to a root in $\textup{supp}(\alpha)$ by an edge in the 
Dynkin diagram for $\gp$.   

\begin{lemma} \label{hess form} 
Let $\alpha \in \Phi^-$.
If $H_{\alpha}^+ = \left\langle E_{\beta}: \beta \in \Phi^+,  \textup{supp}(\beta) \cap 
N(\textup{supp}(\alpha)) \neq 
 \emptyset \right\rangle$, 
$H_{\alpha}^- = \left\langle E_{\beta} : \beta \in \Phi^- \textup{ has } \beta \succeq \alpha \right\rangle$,
and 
$T_{\alpha} = \left\langle [E_{\alpha_i}, E_{-\alpha_i}]: \alpha_i \in
		\textup{supp}(\alpha) \right\rangle$, then
\[H_{\alpha} = H_{\alpha}^- \oplus H_{\alpha}^+ \oplus T_{\alpha}.\]
\end{lemma}

\begin{proof}
Recall that $[E_{\beta}, E_{\gamma}]$ is a nonzero multiple of 
$E_{\beta+\gamma}$ if $\beta+\gamma$ is a root, 
an element $T_{\beta}$ of the Cartan subalgebra if $\gamma = -\beta$,
and zero otherwise.  

This identity implies that 
\[\displaystyle H_{\alpha}^-= \bigcap_{\scriptsize \begin{array}{c}
\textup{Hess. spaces }H \\ \textup{s.t. }E_{\alpha} \in H \end{array}} H \cap \mathfrak{n}^-\] 
and that the Cartan subalgebra $\tor$ intersects
$[H_{\alpha}^-, \bo]$ exactly in $T_{\alpha}$.  The $[\bo,\cdot]$-closure of $T_{\alpha}$
 is $H_{\alpha}^+$.  
% Use Proposition 8.3 of Humphreys Lie Algebras, as well as Jacobi identity (and induct).

We must show that $[H_{\alpha}^-, \bo] \cap \bo \subseteq T_{\alpha} \oplus H_{\alpha}^+$.
Suppose $E_{\gamma} \in \bo$ and $E_{\beta} \in H_{\alpha}^-$ satisfy $\gamma + \beta \in \Phi^+$.
We will find a simple root $\alpha_i \in \textup{supp}(\gamma + \beta) \cap N(\textup{supp}(\alpha))$.
If the support of $\beta$ is contained in the support of $\gamma+\beta$, then any $\alpha_i \in 
\textup{supp}(\beta)$ is as desired, since $\textup{supp}(\beta) \subseteq \textup{supp}(\alpha)$
by definition of $H_{\alpha}^-$.  If $\textup{supp}(\beta) \not \subseteq \textup{supp}(\gamma + \beta)$
then we may write the support of $\gamma$ as the (not necessarily disjoint) union
$\textup{supp}(\gamma) = \textup{supp}(\beta) \cup \textup{supp}(\gamma+\beta)$.  Suppose for every 
$\alpha_i \in \textup{supp}(\gamma+\beta)$ and every $\alpha_j \in \textup{supp}(\alpha)$ we have
$\langle \alpha_i, \alpha_j \rangle = 0$.  Then the support of $\gamma+\beta$ is not connected
to the support of $\alpha$ in the Dynkin diagram, and consequently the support of $\gamma+\beta$
is not connected to the support of $\beta$.  In other words, the support
of $\gamma$ is a disconnected subset of the Dynkin diagram.  This contradicts the fact that
the support of each root is a connected subset of the Dynkin diagram (see \cite[page 169]{B}).
So $E_{\gamma+\beta} \in H_{\alpha}^+$.  \end{proof}

The next two results generalize Lemma \ref{main schubert lemma}.  

\begin{lemma}
If $\alpha$ is a root of the same length as $\theta$, there is a unique maximal Weyl group element $w$ 
that satisfies $w^{-1} \theta = \alpha$.
\end{lemma}

\begin{proof}
Denote the stabilizer of $\theta$ in $W$ by $\St(\theta)$.  Consider the left cosets $\St(\theta) 
\backslash W$.  Each coset has a unique maximal element since
$\St(\theta)$ is parabolic.  Also, each coset $\St(\theta) u$ is determined by the root
$u^{-1} \theta$.  Since $\alpha$ and $\theta$ are in the same $W$-orbit, there exists
$w$ with $w^{-1} \theta = \alpha$.
\end{proof}

\begin{proposition}
Let $\alpha$ be a root of the same length as $\theta$ and let $w$ be the maximal 
Weyl group element with $w^{-1}\theta = \alpha$.
Then $\He(E_{\theta}, H_{\alpha}) = \Cy_w$.
\end{proposition}

\begin{proof}
For each element $u$ in $W$, the flag $[u]$ is in $\He(E_{\theta},H_{\alpha})$ if and only if 
$u^{-1}E_{\theta}u \in H_{\alpha}$.  Since $u^{-1}E_{\theta}u = E_{u^{-1}\theta}$, 
the Hessenberg variety $\He(E_{\theta},H_{\alpha})$ is a union of Schubert cells indexed by the elements in cosets of $\St(\theta) \backslash W$.  We must show that if $E_{u^{-1} \theta} \in
H_{\alpha}$ then $w \geq u$ in the Bruhat order. 

The roots $u^{-1} \theta$ and $w^{-1}\theta$ have the same length.
If $u^{-1}\theta$ and $w^{-1}\theta$ have the same sign, then $u^{-1}\theta \succeq w^{-1}\theta$
if and only if $w \geq u$ by \cite[Proposition 3.2]{St}.   
Now suppose $u^{-1}\theta$ is positive
and $w^{-1}\theta$ is negative.  Without loss of generality, let $u^{-1} \theta=\alpha_i$ 
be simple.  
If $\alpha_i$ is in the support of $w^{-1}\theta$ then $s_i u^{-1}\theta \succeq w^{-1}\theta$ and so
$w \geq us_i$.  Moreover, we know $us_i > u$ since $us_i \alpha_i \in \Phi^-$.  This gives $w \geq u$.

If $\alpha_i$ is not in the support of $w^{-1} \theta$ then there exists an
$\alpha_j \in \textup{supp}(w^{-1}\theta)$ such that $\langle \alpha_i, \alpha_j \rangle \neq 0$
by Lemma \ref{hess form}.   At least one simple root in $\textup{supp}(w^{-1}\theta)$ has
the same length as $w^{-1}\theta$ and hence as $\alpha_i$.  If the Dynkin diagram
for $\gp$ has a multiedge, then the simple roots are long on one side of the multiedge
and short on the other.  So the edge from $\alpha_i$ to $\alpha_j$ cannot be
a multiedge.
This means that $s_j \alpha_i = s_i \alpha_j = \alpha_i+\alpha_j$.
The root $s_ju^{-1}\theta = \alpha_i+\alpha_j$ and so $u > us_j$ 
by \cite[Proposition 3.2]{St}.  Let $u = vs_j$ be a reduced factorization.  Then
$s_js_iv^{-1} \theta = -\alpha_j$ and $w \geq vs_is_j$, again 
by \cite[Proposition 3.2]{St}.  The factorization $vs_is_j$ is reduced because
$vs_i > v$ (\cite[Proposition 3.2]{St}) and $vs_is_j > vs_i$ (because $vs_is_j \alpha_j \in \Phi^-$).
We conclude that $w \geq vs_is_j > u$.
\end{proof}

\section{Questions} \label{questions}

We first ask about the relation between Schubert and Hessenberg varieties.

\begin{question}
Are all Schubert varieties Hessenberg varieties?  If not, describe explicitly the Schubert varieties that are also Hessenberg varieties.
\end{question}

The matrices in the Hessenberg spaces of Corollary \ref{banded form} are 
said to be in {\em banded Hessenberg form},
a form used in numerical analysis (see \cite{dMPS}).  We ask if this algebraic property is related
to the geometric condition of purity.

\begin{question}
Let $X$ be any linear operator. 
If $H$ is in banded Hessenberg form, is the Hessenberg variety $\He(X,H)$ necessarily
pure-dimensional?
\end{question}

The next question arises because of the representations on the cohomology of Springer fibers.   We wonder whether the highest-weight Hessenberg varieties of Sections \ref{highweight} and \ref{generalizations} also carry interesting geometric actions.

\begin{question}
 Does the cohomology of the highest-weight Hessenberg varieties carry interesting group actions?  Is there an interesting group action on the highest-weight Hessenberg variety that permutes its irreducible components?
\end{question}

 \cite{dMPS} proved that regular semisimple Hessenberg varieties are smooth.  
 
\begin{question}
Are all semisimple Hessenberg varieties smooth?
\end{question}

\end{document}